\documentclass[10pt]{article}

\usepackage{amsmath,amsthm,amsfonts,amssymb}
\usepackage{psfig}

\title{Infinite words and universal free actions}
\author{\textsf{Olga Kharlampovich}
\and \textsf{Alexei Myasnikov}
\and \textsf{Denis Serbin}}

\newtheorem{example}{Example}
\newtheorem{corollary}{Corollary}
\newtheorem{prop}{Proposition}
\newtheorem{theorem}{Theorem}
\newtheorem{lemma}{Lemma}

\begin{document}
\maketitle

\begin{abstract}
This is the second paper in a series of three, where we take on the unified theory of non-Archimedean
group actions, length functions and infinite words. Here, for an arbitrary group $G$ of infinite
words over an ordered abelian group $\Lambda$ we construct a $\Lambda$-tree $\Gamma_G$ equipped with
a free action of $G$. Moreover, we show that $\Gamma_G$ is a universal tree for $G$ in the sense that
it isometrically embeds in every $\Lambda$-tree equipped with a free $G$-action compatible with
the original length function on $G$.
\end{abstract}

\section{Introduction}
\label{sec:intro}

The theory of group actions on $\Lambda$-trees goes back to early 1960's. Lyndon introduced abstract 
length functions on groups \cite{Lyndon:1963}, axiomatizing Nielsen cancellation method; he initiated 
the study of groups with real valued length functions. Chiswell related such length functions with 
group actions on $\mathbb{Z}-$ and   $\mathbb{R}$-trees, providing a construction of the tree on 
which the group acts. Tits gave the first formal definition of an $\mathbb{R}$-tree \cite{Tits:1977}. 
In his seminal book \cite{Serre:1980} Serre laid down fundamentals of the theory of groups acting 
freely on simplicial trees. In the following decade Serre's novel approach unified several geometric, 
algebraic, and combinatorial methods of group theory into a unique powerful tool, known today as 
Bass-Serre theory. In their very influential paper \cite{MoSh} Morgan and Shalen linked group 
actions on $\mathbb{R}$-trees with topology and generalized parts of Thurston's Geometrization 
Theorem; they introduced  $\Lambda$-trees for an arbitrary ordered abelian group $\Lambda$ and the 
general form of Chiswell's construction. Thus, it became clear that abstract length functions with 
values in $\Lambda$ and group actions on $\Lambda$-trees are just two equivalent approaches to the 
same realm of group theory questions. The  unified theory was further developed in an important paper 
by Alperin and Bass \cite{AB} where authors state a fundamental problem in the theory of group 
actions on $\Lambda$-trees: find the group theoretic information carried by a $\Lambda$-tree action 
(analogous to Bass-Serre theory), in particular, describe finitely generated groups acting freely 
on $\Lambda$-trees ($\Lambda$-free groups). One of the main breakthroughs in this direction is Rips' 
Theorem, that describes finitely generated $\mathbb{R}$-free groups (see \cite{GLP,BF}). The structure
of finitely generated $\mathbb{Z}^n$-free and $\mathbb{R}^n$-free groups was clarified in 
\cite{Bass:1991,Guirardel:2004}.

Introduction of infinite $\Lambda$-words was one of the major recent developments in the theory of
group actions. In \cite{Myasnikov_Remeslennikov_Serbin:2005} Myasnikov, Remeslennikov and Serbin
showed that groups admitting faithful representations by $\Lambda$-words act freely on some
$\Lambda$-trees, while Chiswell proved the converse \cite{Chiswell:2005}. This gives another equivalent
approach to group actions. Now one can bypass the axiomatic view-point on length functions and work
instead with $\Lambda$-words in the same manner as with ordinary words in standard free 
groups. This allows one to bring into the play naturally and in one package 
such powerful techniques as Nielsen's method, Stallings' graph approach to subgroups, and 
Makanin-Razborov type of elimination processes (see papers 
\cite{Myasnikov_Remeslennikov_Serbin:2005, Myasnikov_Remeslennikov_Serbin:2006, 
Kharlampovich_Myasnikov:2005(2), Kharlampovich_Myasnikov:2006, Kharlampovich_Myasnikov:2010, KMRS2, 
DM, KMS1, KMS2, Nikolaev_Serbin:2011, Nikolaev_Serbin:2011(2), Serbin_Ushakov:2011(1)}). In the case 
when $\Lambda$ is equal to either $\mathbb{Z}^n$ or $\mathbb{Z}^\infty$ all these techniques are 
effective, so many algorithmic problems for $\mathbb{Z}^n$-free groups become decidable, in particular, 
the subgroup membership problem.

In this paper for an arbitrary group $G$ of infinite words over an ordered abelian group $\Lambda$ 
we construct a $\Lambda$-tree $\Gamma_G$ equipped with a free action of $G$. Moreover, we show that 
$\Gamma_G$ is a universal tree for $G$ in the sense that it isometrically embeds in every $\Lambda$-tree 
equipped with a free $G$-action compatible with the original length function on $G$. The construction 
is extremely simple and natural: one just folds every pair of infinite words in $G$ along their 
common initial segments to get the tree $\Gamma_G$. Furthermore, in the case $\Lambda = \mathbb{Z}^n$ 
the construction is effective. Besides, it sheds some light on the nature of the initial Chiswell's 
argument, why it worked and where it came from.

\section{Preliminaries}
\label{sec:prelim}

In this section we introduce basic notions in the theory of $\Lambda$-trees.

\subsection{$\Lambda$-trees}
\label{subs:lambda_trees}

A set $\Lambda$ equipped with addition $+$ and a partial order $\leqslant$ is called a {\em
partially ordered} abelian group if
\begin{enumerate}
\item[(1)] $\langle \Lambda, + \rangle$ is an abelian group,
\item[(2)] $\langle \Lambda, \leqslant \rangle$ is a partially ordered set,
\item[(3)] for all $a,b,c \in \Lambda$, $a \leqslant b$ implies $a + c \leqslant b + c$.
\end{enumerate}

An abelian group $\Lambda$ is called {\em orderable} if there exists a linear order $\leqslant$ on
$\Lambda$, satisfying the condition (3) above. In general, the ordering on $\Lambda$ is not unique.

An ordered abelian group $\Lambda$ is called {\it discretely ordered} if $\Lambda^+$ has a minimal
non-trivial element (we denote it $1_\Lambda$). In this event, for any $a \in \Lambda$ we have
$$a + 1_\Lambda = \min\{b \mid b > a\},\ \ \ a - 1_\Lambda = \max\{b \mid b < a\}.$$

For elements $a,b$ of an ordered group $\Lambda$ the {\it closed segment} $[a,b]$ is defined by
$$[a,b] = \{c \in A \mid a \leqslant c \leqslant b \}.$$

Let $X$ be a non-empty set and $\Lambda$ an ordered abelian group. A {\em $\Lambda$-metric on $X$}
is a mapping $p: X \times X \longrightarrow \Lambda$ such that for all $x,y,z \in X$
\begin{enumerate}
\item[(M1)] $p(x,y) \geqslant 0$,
\item[(M2)] $p(x,y) = 0$ if and only if $x = y$,
\item[(M3)] $p(x,y) = p(y,x)$,
\item[(M4)] $p(x,y) \leqslant p(x,z) + p(y,z)$.
\end{enumerate}

A {\em $\Lambda$-metric space} is a pair $(X,p)$, where $X$ is a non-empty set and $p$ is a
$\Lambda$-metric on $X$. If $(X,p)$ and $(X',p')$ are $\Lambda$-metric spaces, an {\it isometry}
from $(X,p)$ to $(X',p')$ is a mapping $f: X \rightarrow X'$ such that $p(x,y) = p'(f(x),f(y))$ for
all $x,y \in X$.

A {\em segment} in a $\Lambda$-metric space is the image of an isometry $\alpha: [a,b]_\Lambda
\rightarrow X$ for some $a,b \in \Lambda$ and $[a,b]_\Lambda$ is a segment in $\Lambda$. The
endpoints of the segment are $\alpha(a), \alpha(b)$.

We call a $\Lambda$-metric space $(X,p)$ {\em geodesic} if for all $x,y \in X$, there is a segment
in $X$ with endpoints $x,y$ and $(X,p)$ is {\em geodesically linear} if for all $x,y \in X$, there
is a unique segment in $X$ whose set of endpoints is $\{x,y\}$.

A {\em $\Lambda$-tree} is a $\Lambda$-metric space $(X,p)$ such that
\begin{enumerate}
\item[(T1)] $(X,p)$ is geodesic,
\item[(T2)] if two  segments of $(X,p)$ intersect in a single point, which is an endpoint of both,
then their union is a segment,
\item[(T3)] the intersection of two segments with a common endpoint is also a segment.
\end{enumerate}

Let $X$ be a $\Lambda$-tree. We call $e \in X$ an {\em end point} of $X$ if, whenever $e \in [x,y]
\subset X$ either $e = x$ or $e = y$. A {\em linear subtree from $x \in X$} is any linear subtree
$L$ of $X$ having $x$ as an end point. $L$ carries a natural linear ordering with $x$ as least
element. If $y \in L$ then $L_y = \{ z \in L \mid y \leqslant z\}$ is a linear subtree from $y$.

A maximal linear subtree from $x$ is an {\em $X$-ray from $x$}. Observe that (Proposition 2.22
\cite{AB}) that if $L, L'$ are $X$-rays from $x, x'$ respectively such that $L \cap L' \neq
\emptyset$ then $L \cap L'$ is either a closed segment or $L \cap L' = L_v$ for some $v \in X$. In
fact, we call $X$-rays $L$ and $L'$ {\em equivalent} if $L \cap L' = L_v$ for some $v \in X$. The
equivalence classes of $X$-rays for this relation are called {\em ends} of $X$.

We say that group $G$ acts on a $\Lambda$-tree $X$ if every element $g \in G$ defines an isometry
$g : X \rightarrow X$.

Note, that every group has a trivial action on a $\Lambda$-tree, that is, all its elements act as
identity.

Let a group  $G$ act as isometries on a $\Lambda$-tree $X$. $g \in G$ is called {\em elliptic} if
it has a fixed point. $g \in G$ is called an {\em inversion} if it does not have a fixed point, but
$g^2$ does. If $g$ is not elliptic and not an inversion then it is called {\em hyperbolic}. For a
hyperbolic element $g \in G$ define a characteristic set
$$Axis(g) = \{ p \in X \mid [g^{-1} \cdot p, p] \cap [p, g \cdot p] = \{p\} \},$$
which is called the {\em axis of $g$}. $Axis(g)$ meets every $\langle g \rangle$-invariant subtree
of $X$.

A group $G$  acts {\it freely} and {\it without inversions} on a $\Lambda$-tree $X$ if for all $1
\neq g \in G$, $g$ acts as a hyperbolic isometry. In this case we also say that $G$ is
$\Lambda$-free.

\subsection{Length functions}
\label{subs:length}

Let $G$ be a group and $\Lambda$ be an ordered abelian group.  Then a function $l: G \rightarrow
\Lambda$ is called a {\it (Lyndon) length function} on $G$ if the following conditions hold
\begin{enumerate}
\item [(L1)] $\forall\ x \in G:\ l(x) \geqslant 0$ and $l(1) = 0$,
\item [(L2)] $\forall\ x \in G:\ l(x) = l(x^{-1})$,
\item [(L3)] $\forall\ x, y, z \in G:\ c(x,y) > c(x,z) \rightarrow c(x,z) =
c(y,z)$,

\noindent where $c(x,y) = \frac{1}{2}(l(x)+l(y)-l(x^{-1}y))$.
\end{enumerate}

It is not difficult to derive the following two properties of Lyndon length functions from the
axioms (L1)--(L3)
\begin{itemize}
\item $\forall\ x, y \in G:\ l(xy) \leqslant l(x) + l(y)$,
\item $\forall\ x, y \in G:\ 0 \leqslant c(x,y) \leqslant min\{l(x),l(y)\}$.
\end{itemize}
The axiom below helps to describe the connection between $\Lambda$-valued Lyndon length functions
and actions on $\Lambda$-trees.
\begin{enumerate}
\item [(L4)] $\forall\ x \in G:\ c(x,y) \in \Lambda.$
\end{enumerate}

\begin{theorem} \cite{Chiswell:1976}
Let $G$ be a group and $l: G \to \Lambda$ a Lyndon length function satisfying (L4). Then  there
are a $\Lambda$-tree $(X,p)$, an action of $G$ on $X$ and a point $x \in X$ such that $l = l_x$.
\end{theorem}

\subsection{Infinite words}
\label{subs:inf_words}

Let $\Lambda$ be a discretely ordered abelian group with the minimal positive element $1$. It is
going to be clear from the context if we are using $1$ as an element of $\Lambda$, or as an integer.
Let $X = \{x_i \mid i \in I\}$ be a set. Put $X^{-1} = \{x_i^{-1} \mid i \in I\}$ and $X^\pm = X
\cup X^{-1}$. A {\em $\Lambda$-word} is a function of the type
$$w: [1,\alpha_w] \to X^\pm,$$
where $\alpha_w \in \Lambda,\ \alpha_w \geqslant 0$. The element $\alpha_w$ is called the {\em
length} $|w|$ of $w$.

\smallskip

By $W(\Lambda,X)$ we denote the set of all $\Lambda$-words. Observe, that $W(\Lambda,X)$ contains
an empty $\Lambda$-word which we denote by $\varepsilon$.

Concatenation $uv$ of two $\Lambda$-words $u,v \in W(\Lambda,X)$ is an $\Lambda$-word of length
$|u| + |v|$ and such that:
\[ (uv)(a) = \left\{ \begin{array}{ll}
\mbox{$u(a)$}  & \mbox{if $1 \leqslant a \leqslant |u|$} \\
\mbox{$v(a - |u|)$ } & \mbox{if $|u|  < a \leqslant |u| + |v|$}
\end{array}
\right. \]
Next, for any $\Lambda$-word $w$ we define an {\it inverse} $w^{-1}$ as an $\Lambda$-word of the
length $|w|$ and such that
$$w^{-1}(\beta) = w(|w| + 1 - \beta)^{-1} \ \ (\beta \in [1,|w|]).$$

A $\Lambda$-word $w$ is {\it reduced} if $w(\beta + 1) \neq w(\beta)^{-1}$ for each $1 \leqslant
\beta < |w|$. We denote by $R(\Lambda,X)$ the set of all reduced $\Lambda$-words. Clearly,
$\varepsilon \in R(\Lambda,X)$. If the concatenation $uv$ of two reduced $\Lambda$-words $u$ and
$v$ is also reduced then we write $uv = u \circ v$.

\smallskip

For $u \in W(\Lambda,X)$ and $\beta \in [1, \alpha_u]$ by $u_\beta$ we denote the restriction of
$u$ on $[1,\beta]$. If $u \in R(\Lambda,X)$ and $\beta \in [1, \alpha_u]$ then
$$u = u_\beta \circ {\tilde u}_\beta,$$
for some uniquely defined ${\tilde u}_\beta$.

An element ${\rm com}(u,v) \in R(\Lambda,X)$ is called the ({\emph{longest}) {\it common initial
segment} of $\Lambda$-words $u$ and $v$ if
$$u = {\rm com}(u,v) \circ \tilde{u}, \ \ v = {\rm com}(u,v) \circ \tilde{v}$$
for some (uniquely defined) $\Lambda$-words $\tilde{u}, \tilde{v}$ such that $\tilde{u}(1) \neq
\tilde{v}(1)$.

Now, we can define the product of two $\Lambda$-words. Let $u,v \in R(\Lambda,X)$. If ${\rm com}
(u^{-1}, v)$ is defined then
$$u^{-1} = {\rm com}(u^{-1},v) \circ {\tilde u}, \ \ v = {\rm com} (u^{-1},v) \circ {\tilde v},$$
for some uniquely defined ${\tilde u}$ and ${\tilde v}$. In this event put
$$u \ast v = {\tilde u}^{-1} \circ {\tilde v}.$$
The  product ${\ast}$ is a partial binary operation on $R(\Lambda,X)$.

\smallskip

An element $v \in R(\Lambda,X)$ is termed {\it cyclically reduced} if $v(1)^{-1} \neq v(|v|)$. We
say that an element $v \in R(\Lambda,X)$ admits a {\it cyclic decomposition} if $v = c^{-1} \circ u
\circ c$, where $c, u \in R(\Lambda,X)$ and $u$ is cyclically reduced. Observe that a cyclic
decomposition is unique (whenever it exists). We denote by $CR(\Lambda,X)$ the set of all cyclically
reduced words in $R(\Lambda,X)$ and by $CDR(\Lambda,X)$ the set of all words from $R(\Lambda,X)$
which admit a cyclic decomposition.

\smallskip

Below we refer to $\Lambda$-words as {\it infinite words} usually omitting $\Lambda$ whenever it
does not produce any ambiguity.

The following result establishes the connection between infinite words and length functions.
\begin{theorem}
\label{co:3.1} \cite{Myasnikov_Remeslennikov_Serbin:2005}
Let $\Lambda$ be a discretely ordered abelian group and $X$ be a set. Then any subgroup $G$ of
$CDR(\Lambda,X)$ has a free Lyndon length function with values in $\Lambda$ -- the restriction
$L|_G$ on $G$ of the standard length function $L$ on $CDR(\Lambda,X)$.
\end{theorem}

The converse of the theorem above was obtained by I. Chiswell \cite{Chiswell:2005}.

\begin{theorem}
\label{chis} \cite{Chiswell:2005}
Let $G$ have a free Lyndon length function $L : G \rightarrow A$, where $\Lambda$ is a discretely
ordered abelian group. Then there exists a set $X$ and a length preserving embedding $\phi : G
\rightarrow CDR(\Lambda,X)$, that is, $|\phi(g)| = L(g)$ for any $g \in G$.
\end{theorem}

\begin{corollary}
\label{chis-cor} \cite{Chiswell:2005}
Let $G$ have a free Lyndon length function $L : G \rightarrow \Lambda$, where $\Lambda$ is an
arbitrary ordered abelian group. Then there exists an embedding $\phi : G \to CDR(\Lambda',X)$,
where $\Lambda' = \mathbb{Z} \oplus \Lambda$ is discretely ordered with respect to the right
lexicographic order and $X$ is some set, such that, $|\phi(g)| = (0,L(g))$ for any $g \in G$.
\end{corollary}

Theorem \ref{co:3.1}, Theorem \ref{chis}, and Corollary \ref{chis-cor} show that a group has a free
Lyndon length function if and only if it embeds into a set of infinite words and this embedding
preserves the length. Moreover, it is not hard to show that this embedding also preserves regularity
of the length function.

\begin{theorem}
\label{chis-cor-1} \cite{Khan_Myasnikov_Serbin:2007}
Let $G$ have a free regular Lyndon length function $L : G \rightarrow \Lambda$, where $\Lambda$ is
an arbitrary ordered abelian group. Then there exists an embedding $\phi : G \rightarrow R(\Lambda',
X)$, where $\Lambda'$ is a discretely ordered abelian group and $X$ is some set, such that, the
Lyndon length function on $\phi(G)$ induced from $R(\Lambda',X)$ is regular.
\end{theorem}

\section{Universal trees}
\label{sec:universal}

Let $G$ be a subgroup of $CDR(\Lambda,X)$ for some discretely ordered abelian group $\Lambda$ and a
set $X$. We assume $G,\ \Lambda$, and $X$ to be fixed for the rest of this section.

Every element $g \in G$ is a function
$$g: [1,|g|] \rightarrow X^{\pm},$$
with the domain $[1,|g|]$ which a closed segment in $\Lambda$. Since $\Lambda$ can be viewed as a
$\Lambda$-metric space then $[1,|g|]$ is a geodesic connecting $1$ and $|g|$, and every
$\alpha \in [1,|g|]$ we view as a pair $(\alpha, g)$. We would like to identify
initial subsegments of the geodesics corresponding to all elements of $G$ as follows.

\smallskip

Let
$$S_G = \{(\alpha,g) \mid g \in G, \alpha \in [0,|g|]\}.$$
Since for every $f,g \in G$ the word $com(f,g)$ is defined, we can introduce an equivalence relation
on $S_G$ as follows: $(\alpha,f) \sim (\beta,g)$ if and only if $\alpha = \beta \in [0, c(f,g)]$.
Obviously, it is symmetric and reflexive. For transitivity observe that if $(\alpha,f) \sim (\beta,g)$ and
$(\beta,g) \sim (\gamma,h)$ then $0 \leqslant \alpha = \beta = \gamma \leqslant c(f,g), c(g,h)$. Since $c(f,h)
\geqslant \min\{c(f,g), c(g,h)\}$ then $\alpha = \gamma \leqslant c(f,h)$.

Let $\Gamma_G = S_G / \sim$ and $\epsilon = \langle 0, 1 \rangle$, where $\langle \alpha, f \rangle$
is the equivalence class of $(\alpha,f)$.

\begin{prop}
\label{pr:lambda_tree}
$\Gamma_G$ is a $\Lambda$-tree,
\end{prop}
\begin{proof} At first we show that $\Gamma_G$ is a $\Lambda$-metric space. Define the metric by
$$d(\langle \alpha, f \rangle, \langle \beta, g \rangle) = \alpha + \beta - 2\min\{\alpha, \beta, 
c(f,g)\}.$$
Let us check if it is well-defined. Indeed, $c(f,g) \in \Lambda$ is defined for every $f,g \in G$. 
Moreover, let $(\alpha, f) \sim (\gamma,u)$ and $(\beta, g) \sim (\delta,v)$, we want to prove
$$d(\langle \alpha, f \rangle, \langle \beta, g \rangle) = d(\langle \gamma, u \rangle, \langle \delta, 
v \rangle)$$
which is equivalent to
$$\min\{\alpha,\beta,c(f,g)\} = \min\{\alpha,\beta,c(u,v)\}$$
since $\alpha = \gamma,\ \beta = \delta$. Consider the following cases.

\begin{enumerate}
\item[(a)] $\min\{\alpha,\beta\} \leqslant c(u,v)$

Hence, $\min\{\alpha,\beta,c(u,v)\} = \min\{\alpha,\beta\}$ and it is enough to prove
$\min\{\alpha,\beta\}$ $= \min\{\alpha,\beta,c(f,g)\}$. From length function axioms for $G$ we have
$$c(f,g) \geqslant \min\{c(u,f),c(u,g)\},\ \ c(u,g) \geqslant \min\{c(u,v),c(v,g)\}.$$
Hence,
$$c(f,g) \geqslant \min\{c(u,f),c(u,g)\} \geqslant  \min\{c(u,f), \min\{c(u,v),c(v,g)\} \}$$
$$ = \min\{c(u,f),c(u,v),c(v,g)\}.$$
Now, from $(\alpha, f) \sim (\gamma,u),\ (\beta, g) \sim (\delta,v)$ it follows that
$\alpha \leqslant c(u,f),\ \beta \leqslant c(v,g)$ and combining it with the assumption
$\min\{\alpha,\beta\} \leqslant c(u,v)$ we have
$$c(f,g) \geqslant \min\{c(u,f),c(u,v),c(v,g)\} \geqslant \min\{\alpha,\beta\},$$
or, in other words,
$$\min\{\alpha,\beta,c(f,g)\} = \min\{\alpha,\beta\}.$$

\item[(b)] $\min\{\alpha,\beta\} > c(u,v)$

Hence, $\min\{\alpha,\beta,c(u,v)\} = c(u,v)$ and it is enough to prove $c(f,g) = c(u,v)$.

Since
$$c(u,f) \geqslant \alpha > c(u,v),\ c(v,g) \geqslant \beta > c(u,v),$$
then $\min\{c(u,f),c(u,v),c(v,g)\} = c(u,v)$ and
$$c(f,g) \geqslant \min\{c(u,f),c(u,v),c(v,g)\} = c(u,v).$$
Now we prove that $c(f,g) \leqslant c(u,v)$. From length function axioms for $G$ we have
$$c(u,v) \geqslant \min\{c(v,g),c(u,g)\} = c(u,g) \geqslant \min\{c(v,g),c(u,v)\} = c(u,v),$$
that is, $c(u,v) = c(u,g)$. Now,
$$c(u,v) = c(u,g) \geqslant \min\{c(u,f),c(f,g)\},$$
where $\min\{c(u,f),c(f,g)\} = c(f,g)$ since otherwise we have $c(u,v) \geqslant c(u,f) \geqslant \alpha$ -
a contradiction. Hence, $c(u,v) \geqslant c(f,g)$ and we have $c(f,g) = c(u,v)$.
\end{enumerate}

By definition of $d$, for any $\langle \alpha, f \rangle,\ \langle \beta, g \rangle$
we have
$$d(\langle \alpha, f \rangle, \langle \beta, g \rangle) = d(\langle \beta, g \rangle, \langle 
\alpha, f \rangle) \geqslant 0,$$
$$d(\langle \alpha, f \rangle, \langle \alpha, f \rangle) = 0.$$
If
$$d(\langle \alpha, f \rangle, \langle \beta, g \rangle) = \alpha + \beta - 2 \min\{\alpha, \beta, 
c(f,g)\} = 0$$
then $\alpha + \beta = 2 \min\{\alpha,\beta,c(f,g)\}$. It is possible only if $\alpha = \beta \leqslant 
c(f,g)$ which implies $\langle \alpha, f \rangle = \langle \beta, g \rangle$. Finally, we have to 
prove the triangle inequality
$$d(\langle \alpha, f \rangle, \langle \beta, g \rangle) \leqslant d(\langle \alpha, f \rangle, \langle 
\gamma, h \rangle) + d(\langle \beta, g \rangle, \langle \gamma, h \rangle)$$
for every $\langle \alpha, f \rangle,\ \langle \beta, g \rangle,\ \langle \gamma, h \rangle \in 
\Gamma_G$. The inequality above is equivalent to
$$\alpha + \beta - 2\min\{\alpha, \beta, c(f,g)\} \leqslant \alpha + \gamma$$
$$ - 2 \min\{\alpha, \gamma, c(f,h) + \beta + \gamma - 2\min\{\beta,\gamma,c(g,h)\}\}$$
which comes down to
$$\min\{\alpha,\gamma,c(f,h)\} + \min\{\beta,\gamma,c(g,h)\} \leqslant  \min\{\alpha,\beta,c(f,g)\} + 
\gamma.$$

\smallskip

First of all, observe that for any $\alpha,\beta,\gamma \in \Lambda$ the triple $(\min\{\alpha, 
\beta\},$ $\min\{\alpha,\gamma\},$ $\min\{\beta,\gamma\})$ is isosceles. Hence, by Lemma 1.2.7(1) 
\cite{Chiswell:2001}, the triple
$$(\min\{\alpha,\beta,c(f,g)\},\ \min\{\alpha,\gamma,c(f,h)\},\ \min\{\beta,\gamma,c(g,h)\})$$
is isosceles too. In particular,
$$\min\{\alpha,\beta,c(f,g)\} \geqslant \min\{ \min\{\alpha,\gamma,c(f,h)\},\ \min\{\beta,\gamma,c(g,h)\} \}$$
$$ = \min\{\alpha,\beta,\gamma,c(f,h),c(g,h)\}.$$
Now, if
$$\min\{\alpha,\beta,\gamma,c(f,h),c(g,h)\} = \min\{\alpha,\gamma,c(f,h)\}$$
then $\min\{\beta,\gamma,c(g,h)\} = \gamma$ and
$$\min\{\alpha,\gamma,c(f,h)\} + \min\{\beta,\gamma,c(g,h)\} \leqslant  \min\{\alpha,\beta,c(f,g)\} + \gamma$$
holds. If
$$\min\{\alpha,\beta,\gamma,c(f,h),c(g,h)\} = \min\{\beta,\gamma,c(g,h)\}$$
then $\min\{\alpha,\gamma,c(f,h)\} = \gamma$ and
$$\min\{\alpha,\gamma,c(f,h)\} + \min\{\beta,\gamma,c(g,h)\} \leqslant  \min\{\alpha,\beta,c(f,g)\} + \gamma$$
holds again. So, $d$ is a $\Lambda$-metric.

\smallskip

Finally, we want to prove that $\Gamma_G$ is $0$-hyperbolic with respect to $\epsilon = \langle 0, 
1 \rangle$ (and, hence, with respect to any other point in $\Gamma_G$). It is enough to prove that 
the triple
$$((\langle \alpha, f \rangle \cdot \langle \beta, g \rangle)_\epsilon,\ (\langle \alpha, f \rangle 
\cdot \langle \gamma, h \rangle)_\epsilon,\ (\langle \beta, g \rangle \cdot \langle \gamma, h 
\rangle)_\epsilon)$$
is isosceles for every $\langle \alpha, f \rangle,\ \langle \beta, g \rangle,\ \langle \gamma, h 
\rangle \in \Gamma_G$. But by definition of $d$ the above triple is isosceles if and only if
$$(\min\{\alpha,\beta,c(f,g)\},\ \min\{ \alpha, \gamma, c(f,h) \},\ \min\{ \beta, \gamma, c(g,h) \})$$
is isosceles which holds.

\smallskip

So, $\Gamma_G$ is a $\Lambda$-tree.

\end{proof}

Since $G$ is a subset of $CDR(\Lambda,X)$ and every element $g \in G$ is a function defined on 
$[1_A,|g|]$ with values in $X^\pm$ then we can define a function
$$\xi : (\Gamma_G - \{\epsilon\}) \rightarrow X^\pm,\ \ \xi(\langle \alpha, g \rangle) = g(\alpha).$$
It is easy to see that $\xi$ is well-defined. Indeed, if $(\alpha, g) \sim (\alpha_1,g_1)$
then $\alpha = \alpha_1 \leqslant c(g,g_1)$, so $g(\alpha) = g_1(\alpha_1)$. Moreover, since every $g \in G$
is reduced then $\xi(p) \neq \xi(q)^{-1}$ whenever $d(p,q) = 1_A$.

$\xi$ can be extended to a function
$$\Xi : geod(\Gamma_G)_\epsilon \to R(\Lambda,X),$$
where $geod(\Gamma_G)_\epsilon = \{ (\epsilon,p] \mid p \in \Gamma_G \}$, so that
$$\Xi(\ (\epsilon, \langle \alpha, g \rangle]\ )(t) = g(t),\ t \in [1_A,\alpha].$$
That is, $\Xi(\ (\epsilon, \langle \alpha, g \rangle]\ )$ is the initial subword of $g$ of length 
$\alpha$, and
$$\Xi(\ (\epsilon, \langle |g|, g \rangle]\ ) = g.$$
On the other hand, if $g \in G$ and $\alpha \in [1_A,|g|]$ then
the initial subword of $g$ of length $\alpha$ uniquely corresponds to $\Xi(\ (\epsilon, \langle 
\alpha, g \rangle]\ )$. If $(\alpha, g) \sim (\alpha_1,g_1)$ then $\alpha = \alpha_1 \leqslant c(g,g_1)$, 
and since $g(t) = g_1(t)$ for any $t \in [1_A,c(g,g_1)]$ then
$$\Xi(\ (\epsilon, \langle \alpha, g \rangle]\ ) = \Xi(\ (\epsilon, \langle \alpha_1, g_1 \rangle]\ ).$$

\begin{lemma}
\label{le:subword_prod}
Let $u,v \in R(\Lambda,X)$. If $u \ast v$ is defined then $u \ast a$ is also defined, where $v = a 
\circ b$. Moreover, $u \ast a$ is an initial subword of either $u$ or $u \ast v$.
\end{lemma}
\begin{proof} The proof follows from Figure \ref{pic1}.

\begin{figure}[htbp]
\label{pic1}
\centering{\mbox{\psfig{figure=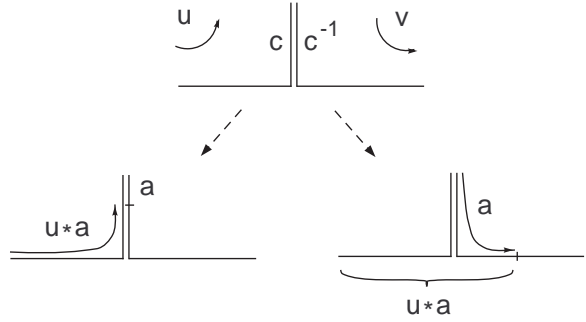,height=2in}}}
\caption{Possible cancellation diagrams in Lemma \ref{le:subword_prod}.}
\end{figure}

\end{proof}

Now, since for every $\langle \alpha, g \rangle \in \Gamma_G,\ \Xi(\ (\epsilon, \langle \alpha, g 
\rangle]\ )$ is an initial subword of $g \in G$ then by Lemma \ref{le:subword_prod}, $f \ast 
\Xi(\ (\epsilon, \langle \alpha, g \rangle]\ )$ is defined for any $f \in G$. Moreover, again by 
Lemma \ref{le:subword_prod},  $f \ast \Xi(\ (\epsilon, \langle \alpha, g \rangle]\ )$ is an initial 
subword of either $f$ or $f \ast g$. More precisely,
$$f \ast \Xi(\ (\epsilon, \langle \alpha, g \rangle]\ ) = \Xi(\ (\epsilon, \langle |f| - \alpha, 
f \rangle]\ )$$
if $f \ast \Xi(\ (\epsilon, \langle \alpha, g \rangle]\ )$ is an initial subword of $f$, and
$$f \ast \Xi(\ (\epsilon, \langle \alpha, g \rangle]\ ) = \Xi(\ (\epsilon, \langle |f| + \alpha - 
2 c(f^{-1},g), f\ast g \rangle]\ )$$
if $f \ast \Xi(\ (\epsilon, \langle \alpha, g \rangle]\ )$ is an initial subword of $f \ast g$.

\smallskip

Hence, we define a (left) action of $G$ on $\Gamma_G$ as follows:
$$f \cdot \langle \alpha, g \rangle = \langle |f| + \alpha - 2\min\{\alpha,c(f^{-1},g)\}, f \rangle$$
if $\alpha \leqslant c(f^{-1},g)$, and
$$f \cdot \langle \alpha, g \rangle = \langle |f| + \alpha - 2\min\{\alpha,c(f^{-1},g)\}, f\ast g 
\rangle$$
if $\alpha > c(f^{-1},g)$.

\smallskip

The action is well-defined. Indeed, it is easy to see that $f \cdot \langle \alpha, g \rangle = f 
\cdot \langle \alpha_1, g_1 \rangle$ whenever $(\alpha, g) \sim (\alpha_1,g_1)$.

\begin{lemma}
\label{le:isometric}
The action of $G$ on $\Gamma_G$ defined above is isometric.
\end{lemma}
\begin{proof} Observe that it is enough to prove
$$d(\epsilon, \langle \alpha, g \rangle) = d(f \cdot \epsilon, f \cdot \langle \alpha, g \rangle)$$
for every $f,g \in G$. Indeed, from the statement above it is going to follow that the geodesic tripod
$(\epsilon,\langle |g|,g \rangle, \langle |h|,h \rangle)$ is isometrically mapped to the geodesic tripod
$(\langle |f|,f \rangle, f \cdot \langle |g|,g \rangle, f \cdot \langle |h|,h \rangle)$ and isometricity
follows.

We have
$$d(\epsilon, \langle \alpha, g \rangle) = d(\langle 0, 1 \rangle, \langle \alpha, g \rangle) = 0 + 
\alpha - 2 \min\{ 0, \alpha, c(1,g) \} = \alpha,$$
$$d(f \cdot \epsilon, f \cdot \langle \alpha, g \rangle) = d(\langle |f|,f \rangle, f \cdot \langle 
\alpha, g \rangle).$$
Consider two cases.

\begin{enumerate}
\item[(a)] $\alpha \leqslant c(f^{-1},g)$

Hence,
$$d(\langle |f|,f \rangle, f \cdot \langle \alpha, g \rangle) = d(\langle |f|,f \rangle, \langle |f| 
- \alpha, f \rangle)$$
$$ = |f| + |f| - \alpha - 2 \min\{|f|,|f|-\alpha,c(f,f)\} = |f| + |f|-\alpha - 2(|f|-\alpha) = \alpha.$$

\item[(b)] $\alpha > c(f^{-1},g)$

Hence,
$$d(\langle |f|,f \rangle, f \cdot \langle \alpha, g \rangle) = d(\langle |f|,f \rangle, \langle |f| 
+ \alpha - 2c(f^{-1},g), f\ast g \rangle)$$
$$ = |f| + |f| + \alpha - 2c(f^{-1},g) - 2\min\{|f|,|f| + \alpha - 2c(f^{-1},g),c(f,f\ast g\})$$
$$= 2|f| + \alpha - 2c(f^{-1},g) - 2\min\{|f| + \alpha - 2c(f^{-1},g),c(f,f\ast g)\}.$$
Let $f = f_1 \circ c^{-1},\ g = c \circ g_1,\ |c| = c(f^{-1},g)$. Then $|f| + \alpha - 2c(f^{-1},g) 
= |f_1|+\alpha-c(f^{-1},g) > |f_1|$. At the same time, $c(f,f \ast g) = |f_1|$, so $\min\{|f| + 
\alpha - 2c(f^{-1},g),c(f,f\ast g)\} = |f_1|$ and
$$d(\langle |f|,f \rangle, f \cdot \langle \alpha, g \rangle) = 2|f| + \alpha - 2c(f^{-1},g) - 2|f_1| 
= 2|f| + \alpha - 2|c| - 2|f_1| = \alpha$$.

\end{enumerate}

\end{proof}

\begin{prop}
\label{pr:action}
The action of $G$ on $\Gamma_G$ defined above is free and $L_\epsilon(g) = |g|$. Moreover, $\Gamma_G$ 
is minimal with respect to this action if and only if $G$ contains a cyclically reduced element $h 
\in G$, that is, $|h^2| = 2|h|$.
\end{prop}
\begin{proof} {\bf Cialm 1.} The stabilizer of every $x \in \Gamma_G$ is trivial.

\smallskip

Next, suppose $f \cdot \langle \alpha, g \rangle = \langle \alpha, g \rangle$.
First of all, if $\alpha = 0$ then $|f| + \alpha - 2\min\{\alpha,c(f^{-1},g)\} = |f|$ then $|f| = 
\alpha = 0$. Also, if $c(f^{-1},g) = 0$ then $|f| + \alpha - 2\min\{\alpha,c(f^{-1},g)\} = |f| + 
\alpha$ which has to be equal to $\alpha$ form our assumption. In both cases $f = 1$ follows.

Assume $f \neq 1$ (which implies $\alpha,\ c(f^{-1},g) \neq 0$) and consider the following cases.

\begin{enumerate}
\item[(a)] $\alpha < c(f^{-1},g)$

Hence, from
$$\langle \alpha, g \rangle = \langle |f| - \alpha, f \rangle$$
we get $\alpha = |f| - \alpha \leqslant c(f,g)$. In particular, $|f| = 2 \alpha$.

Consider the product $f \ast g$. We have
$$f = f_1 \circ com(f^{-1},g)^{-1},\ g = com(f^{-1},g) \circ g_1.$$
Since $\alpha < c(f^{-1},g)$ then we have $com(f^{-1},g) = c_\alpha \circ c,\ |c_\alpha| = \alpha$. 
Hence,
$$f = f_1 \circ c^{-1} \circ c_\alpha^{-1},\ g = c_\alpha \circ c \circ g_1.$$
On the other hand, from $|f| = 2\alpha$ we get $|f_1| + |c| = \alpha \leqslant c(f,g)$, so, $com(f,g)$ 
has $f_1 \circ c$ as initial subword. That is, $g = f_1 \circ c \circ g_2$, but now comparing two 
representations of $g$ above we get $c_\alpha = f_1 \circ c^{-1}$ and $c_\alpha \ast c \neq c_\alpha 
\circ c$ - a contradiction.

\item[(b)] $\alpha = c(f^{-1},g)$

We have $f = f_1 \circ c_\alpha^{-1},\ g = c_\alpha \circ g_1,\ |c_\alpha| = \alpha$. From  $\langle 
\alpha, g \rangle = \langle |f| - \alpha, f \rangle$ we get $\alpha = |f| - \alpha \leqslant c(f,g)$, so 
$|f| = 2\alpha$ and $|f_1| = \alpha$. Since $|f_1| = \alpha \leqslant c(f,g)$ then $g = f_1 \circ g_2$ 
from which it follows that $f_1 = c_\alpha$. But then $f_1 \ast c_\alpha^{-1} \neq f_1 \circ 
c_\alpha^{-1}$ - contradiction.

\item[(c)] $\alpha > c(f^{-1},g)$

Hence, from
$$\langle \alpha, g \rangle = \langle |f| + \alpha - 2 c(f^{-1},g), f \ast g \rangle$$
we get $\alpha = |f| + \alpha - 2 c(f^{-1},g) \leqslant c(g,f \ast g)$. In particular, $|f| = 2 c(f^{-1}, 
g)$.

Consider the product $f \ast g$. We have
$$f = f_1 \circ c^{-1},\ g = c \circ g_1,$$
where $c = com(f^{-1},g)$. Hence, $|f_1| = |c| < \alpha \leqslant c(g,f \ast g) = c(g, f_1 \circ g_1)$.
It follows that $g = f_1 \circ g_2$ and, hence, $c = f_1$ which is impossible.

\end{enumerate}

\smallskip

{\bf Cialm 2.} $L_\epsilon(g) = |g|$

\smallskip

We have $L_\epsilon(g) = d(\epsilon, g \cdot \epsilon)$. Hence, by definition of $d$
$$d(\langle 0,1\rangle, g \cdot \langle 0,1\rangle) = d(\langle 0, 1 \rangle, \langle |g|, g \rangle) 
= 0 + |g| - 2 \min\{ 0, |g|, c(1,g) \} = |g|.$$

\smallskip

{\bf Cialm 3.} $\Gamma_G$ is minimal with respect to the action if and only if $G$ contains a 
cyclically reduced element $h \in G$, that is, $|h^2| = 2|h|$.

\smallskip

Suppose there exists a cyclically reduced element $h \in G$. Let $\Delta \subset \Gamma_G$ be a 
$G$-invariant subtree.

First of all, observe that $\epsilon \notin \Delta$. Indeed, if $\epsilon \in \Delta$ then $f \cdot 
\epsilon \in \Delta$ for every $f \in G$ and since $\Delta$ is a tree then $[\epsilon, f \cdot 
\epsilon] \in \Delta$ for every $f \in G$. At the same time, $\Gamma_G$ is spanned by $[\epsilon, f 
\cdot \epsilon],\ f \in G$, so, $\Delta = \Gamma_G$ - a contradiction.

Let $u \in \Delta$. By definition of $\Gamma_G$ there exists $g \in G$ such that $u \in [\epsilon, 
g \cdot \epsilon]$. Observe that $A_g \subseteq \Delta$. Indeed, for example by Theorem 1.4 
\cite{Chiswell:2001}, if $[u,p]$ is the bridge between $u$ and $A_g$ then $p = Y(g^{-1} \cdot u,\ 
u,\ g \cdot u)$. In particular, $p \in \Delta$ and since for every $v \in A_g$ there exist $g_1, 
g_2 \in C_G(g)$ such that $v \in [g_1 \cdot p,\ g_2 \cdot p]$ then $A_g \subseteq \Delta$.

Observe that if $g$ is cyclically reduced then $\epsilon \in A_g$, that is, $\epsilon \in \Delta$ -
a contradiction. More generally, $\Delta \cap A_f = \emptyset$ for every cyclically reduced $f \in G$.
Hence, let $[p,q]$ be the bridge between $A_g$ and $A_h$ so that $p \in A_g,\ q \in A_h$. Then by
Lemma 2.2 \cite{Chiswell:2001}, $[p,q] \subset A_{gh}$, in particular, $p,q \in A_{gh}$. It follows 
that $A_{gh} \subseteq \Delta,\ q \in A_{gh} \cap A_h$, and $\Delta \cap A_h \neq \emptyset$ - a 
contradiction.

Hence, there can be no proper $G$-invariant subtree $\Delta$.

\smallskip

Now, suppose $G$ contains no cyclically reduced element. Hence, $\epsilon \notin A_f$ for every $f 
\in G$. Let $\Delta$ be spanned by $A_f,\ f \in G$. Obviously, $\Delta$ is $G$-invariant. Indeed, 
let $u \in [p,q]$, where $p \in A_f,\ q \in A_g$ for some $f,g \in G$. Then $h \cdot u \in [h \cdot p, 
h \cdot q]$, where $h \cdot p \in h \cdot A_f = A_{hfh^{-1}},\ h \cdot q \in h \cdot A_g = A_{hgh^{-1}}$, 
that is, $h \in \Delta$.

Finally, $\epsilon \in \Gamma_G - \Delta$.

\end{proof}

\begin{prop}
\label{pr:universal}
If $(Z,d')$ is a $\Lambda$-tree on which $G$ acts freely as isometries, and $w \in Z$ is such that
$L_w(g) = |g|,\ g \in G$ then there is a unique $G$-equivariant isometry $\mu: \Gamma_G \to Z$
such that $\mu(\epsilon) = w$, whose image is the subtree of $Z$ spanned by the orbit $G \cdot w$
of $w$.
\end{prop}
\begin{proof} Define a mapping $\mu : \Gamma_G \rightarrow Z$ as follows
$$\mu(\langle \alpha, f \rangle) = x\ {\rm if}\ d'(w,x) = \alpha,\ d'(f \cdot w,x) = |f|-\alpha.$$
Observe that $\mu(\epsilon) = \mu(\langle 0,1\rangle) = w$

\smallskip

{\bf Claim 1.} $\mu$ is an isometry.

\smallskip

Let $\langle \alpha, f \rangle,\ \langle \beta, g \rangle \in \Gamma_G$. Then by definition of $d$
we have
$$d(\langle \alpha, f \rangle, \langle \beta, g \rangle) = \alpha + \beta - 2 \min\{ \alpha, \beta, 
c(f,g)\}.$$ 
Let $x = \mu(\langle \alpha, f \rangle),\ y = \mu(\langle \beta, g \rangle)$. Then By Lemma 1.2 
\cite{Chiswell:2001} in $(Z,d')$ we have
$$d'(x,y) = d(w,x) + d(w,y) - 2\min\{d(w,x),d(w,y),d(w,z)\},$$
where $z = Y(w,f\cdot w, g\cdot w)$. Observe that $d(w,x) = \alpha,\ d(w,y) = \beta$. At the same 
time, since $L_w(g) = |g|,\ g \in G$ then
$$d(w,z) = \frac{1}{2}(d(w,f\cdot w) + d(w,g\cdot w) - d(f\cdot w,g\cdot w)) = \frac{1}{2}(|f| + |g| 
- |f^{-1} g|) = c(f,g),$$
and
$$d(\mu(\langle \alpha, f \rangle), \mu(\langle \beta, g \rangle)) = d'(x,y) = \alpha + \beta -
2\min\{\alpha,\beta,c(f,g)\} = d(\langle \alpha, f \rangle, \langle \beta, g \rangle).$$

\smallskip

{\bf Claim 1.} $\mu$ is equivariant.

\smallskip

We have to prove
$$\mu(f\cdot \langle \alpha, g \rangle) = f \cdot \mu(\langle \alpha, g \rangle).$$
Let $x = \mu(\langle \alpha, g \rangle),\ y = \mu(f\cdot \langle \alpha, g \rangle)$.
By definition of $\mu$ we have $d'(w,x) = \alpha,\ d'(g \cdot w,x) = |g|-\alpha$.
\begin{enumerate}
\item[(a)] $\alpha \leqslant c(f^{-1},g)$

Hence,
$$f \cdot \langle \alpha, g \rangle = \langle |f| - \alpha, f \rangle.$$
and to prove $y = f\cdot x$ it is enough to show that $d'(w,f\cdot x) = |f|-\alpha$
and $d'(f\cdot w,f\cdot x) = \alpha$.

Observe that the latter equality holds since $d'(f\cdot w,f\cdot x) = d'(w,x) = \alpha$.
To prove the former one, by Lemma 1.2 \cite{Chiswell:2001} we have
$$d(w,f\cdot x) = d'(w,f\cdot w) + d'(f\cdot x, f\cdot w)$$
$$ - 2\min\{d'(w,f\cdot w),
d'(f\cdot x, f\cdot w), d'(f\cdot w,z)\},$$
where $z = Y(w,f\cdot w, (fg)\cdot w)$. Also,
$$d'(f\cdot w,z) = \frac{1}{2}(d'(f\cdot w,w) + d'(f\cdot w,(fg)\cdot w)- d'(w,(fg)\cdot w))$$
$$ = \frac{1}{2}(|f|+|g|-|f^{-1}g|) = c(f^{-1},g).$$
Since, $d'(w,f\cdot w) = |f|,\ d'(f\cdot x, f\cdot w) = \alpha$ then $\min\{d'(w,f\cdot w),
d'(f\cdot x, f\cdot w), d'(f\cdot w,z)\} = \alpha$, and
$$d'(w,f\cdot x) = |f|+\alpha - 2\alpha = |f|-\alpha.$$

\item[(b)] $\alpha > c(f^{-1},g)$

Hence,
$$f \cdot \langle \alpha, g \rangle = \langle |f| + \alpha - 2c(f^{-1},g), f\ast g \rangle.$$
and to prove $y = f\cdot x$ it is enough to show that $d'(w,f\cdot x) = |f| + \alpha - 2c(f^{-1},g)$
and $d'(f\cdot x,(fg) \cdot w) = |fg|-(|f| + \alpha - 2c(f^{-1},g))$.

Observe that $d'(f\cdot x,(fg) \cdot w) = d'(x,gw) = |g| - \alpha = |fg|-(|f| + \alpha - 2c(f^{-1},g))$,
so the latter equality holds.

By Lemma 1.2 \cite{Chiswell:2001} we have
$$d(w,f\cdot x) = d'(w,f\cdot w) + d'(f\cdot x, f\cdot w)$$
$$ - 2\min\{d'(w,f\cdot w),
d'(f\cdot x, f\cdot w), d'(f\cdot w,z)\},$$
where $z = Y(w,f\cdot w, (fg)\cdot w)$. Also,
$$d'(f\cdot w,z) = \frac{1}{2}(d'(f\cdot w,w) + d'(f\cdot w,(fg)\cdot w)- d'(w,(fg)\cdot w))$$
$$ = \frac{1}{2}(|f|+|g|-|f^{-1}g|) = c(f^{-1},g).$$
$d'(w,f\cdot w) = |f|,\ d'(f\cdot x, f\cdot w) = \alpha$, so $\min\{d'(w,f\cdot w),
d'(f\cdot x, f\cdot w), d'(f\cdot w,z)\} = d'(f\cdot w,z) = c(f^{-1},g)$, and
$$d(w,f\cdot x) = |f| + \alpha - 2c(f^{-1},g).$$

\smallskip

{\bf Claim 1.} $\mu$ is unique.

\smallskip

Observe that if $\mu' : \Gamma_G \rightarrow Z$ is another equivariant isometry such that
$\mu'(\epsilon) = w$ then for every $g \in G$ we have
$$\mu'(\langle |g|,g\rangle) = \mu'(g\cdot \langle 0,1\rangle) = g \cdot \mu'(\langle 0,1\rangle) =
g\cdot w.$$
That is, $\mu'$ agrees with $\mu$ on $G\cdot \epsilon$, hence $\mu = \mu'$ because isometries
preserve geodesic segments.

Thus, $\mu$ is unique. Moreover, $\mu(\Gamma_G)$ is the subtree of $Z$ spanned by $G \cdot w$.

\end{enumerate}

\end{proof}

\section{Examples}
\label{sec:examples}

Here we consider two examples of subgroups of $CDR(\Lambda, X)$, where $\Lambda = \mathbb{Z}^2$ and
$X$ an arbitrary alphabet, and explicitly construct the corresponding universal trees for these groups.

\begin{example}
\label{ex:1}
Let $F = F(X)$ be a free group with basis $X$ and the standard length function $|\cdot|$, and let
$u \in F$ a cyclically reduced element which is not a proper power. If we assume that $\mathbb{Z}^2
= \langle 1, t \rangle$ is the additive group of linear polynomials in $t$ ordered lexicographically
then the HNN-extension
$$G = \langle F, s \mid u^s = u \rangle$$
embeds into $CDR(\mathbb{Z}^2, X)$ under the following map $\phi$:
$$\phi(x) = x,\ \forall\ x \in X,$$
\[ \mbox{$\phi(s)(\beta)$} = \left\{ \begin{array}{ll}
\mbox{$u(\alpha)$,} & \mbox{if $\beta = m |u| + \alpha, m \geqslant 0, 1 \leqslant \alpha \leqslant |u|$,} \\
\mbox{$u(\alpha)$,} & \mbox{if $\beta = t - m |u| + \alpha, m > 0, 1 \leqslant \alpha \leqslant |u|$.}
\end{array}
\right.
\]

\begin{figure}[htbp]
\label{pic2}
\centering{\mbox{\psfig{figure=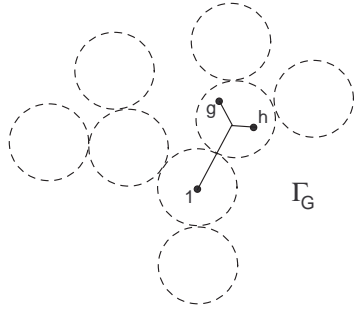,height=1.8in}}}
\caption{$\Gamma_G$ as a $\mathbb{Z}$-tree of $\mathbb{Z}$-trees.}
\end{figure}

It is easy to see that $|\phi(s)| = t$ and $\phi(s)$ commutes with $u$ in $CDR(\mathbb{Z}^2, X)$. To simplify
the notation we identify $G$ with its image $\phi(G)$.

Every element $g$ of $G$ can be represented as the following reduced $\mathbb{Z}^2$-word
$$g = g_1 \circ s^{\delta_1} \circ g_2 \circ \cdots \circ g_k \circ s^{\delta_k} \circ g_{k+1},$$
where $[g_i, u] \neq 1$. Now, according to the construction described in Section \ref{sec:universal}, the
universal tree $\Gamma_G$ consists of the segments in $\mathbb{Z}^2$ labeled by elements from $G$ which
are glued together along their common initial subwords.

Thus, $\Gamma_G$ can be viewed as a $\mathbb{Z}$-tree
of $\mathbb{Z}$-trees which are Cayley graphs of $F(X)$ and every vertex $\mathbb{Z}$-subtree can
be associated with a right representative in $G$ by $F$. The end-points of the segments $[1, |g|]$ and
$[1, |h|]$ labeled respectively by $g$ and $h$ belong to the same vertex $\mathbb{Z}$-subtree if and only if
$h^{-1} g \in F$.

\begin{figure}[htbp]
\label{pic3}
\centering{\mbox{\psfig{figure=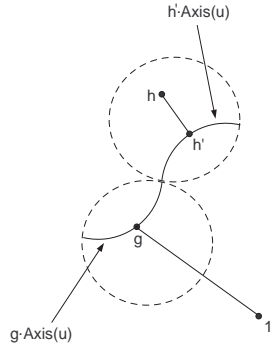,height=2in}}}
\caption{Adjacent $\mathbb{Z}$-subtrees in $\Gamma_G$.}
\end{figure}

In other words, $\Gamma_G$ is a ``more detailed'' version of the Bass-Serre tree $T$ for $G$,
in which every vertex is replaced by the Cayley graph of the base group $F$ and the adjacent $\mathbb{Z}$-subtrees
of $\Gamma_G$ corresponding to the representatives $g$ and $h$ are ``connected'' by means of $s^{\pm}$ which
extends $g \cdot Axis(u)$ to $h' \cdot Axis(u)$, where $h'^{-1} h \in F$ and $g^{-1} h \in s^{\pm} F$.

\end{example}

The following example is a generalization of the previous one.

\begin{example}
\label{ex:2}
Let $F = F(X)$ be a free group with basis $X$ and the standard length function $|\cdot|$, and let
$u, v \in F$ be cyclically reduced elements which is not a proper powers and such that $|u| = |v|$.
The HNN-extension
$$H = \langle F, s \mid u^s = v \rangle$$
embeds into $CDR(\mathbb{Z}^2, X)$ under the following map $\psi$:
$$\psi(x) = x,\ \forall\ x \in X,$$
\[ \mbox{$\psi(s)(\beta)$} = \left\{\begin{array}{ll}
\mbox{$u(\alpha)$,} & \mbox{if $\beta = m |u| + \alpha, m \geqslant 0, 1 \leqslant \alpha \leqslant |u|$,} \\
\mbox{$v(\alpha)$,} & \mbox{if $\beta = t - m |v| + \alpha, m > 0, 1 \leqslant \alpha \leqslant |v|$.}
\end{array}
\right.
\]
It is easy to see that $|\psi(s)| = t$ and $u \circ \psi(s) = \psi(s) \circ v$ in $CDR(\mathbb{Z}^2, X)$.
Again, to simplify the notation we identify $H$ with its image $\psi(H)$.

The structure of $\Gamma_H$ is basically the same as the structure of $\Gamma_G$ in Example \ref{ex:1}.
The only difference is that the adjacent $\mathbb{Z}$-subtrees of $\Gamma_H$ corresponding to the representatives
$g$ and $h$ are ``connected'' by means of $s^{\pm}$ which extends $g \cdot Axis(u)$ to $h' \cdot Axis(v)$, where
$h'^{-1} h \in F$ and $g^{-1} h \in s^{\pm} F$.

\begin{figure}[htbp]
\label{pic4}
\centering{\mbox{\psfig{figure=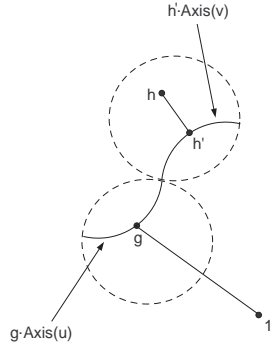,height=2in}}}
\caption{Adjacent $\mathbb{Z}$-subtrees in $\Gamma_H$.}
\end{figure}
\end{example}

\section{Effective $\Lambda$-trees}
\label{se:effect}

In this section we introduce some basic notions concerning effectiveness when dealing
with groups of infinite words and corresponding universal trees.

\subsection{Infinite words viewed as computable functions}
\label{subs:comp_func}

We say that a group $G = \langle Y \rangle,\ Y = \{y_1, \ldots, y_m\}$ has an {\em effective
representation by $\Lambda$-words over an alphabet $X$} if $G \subset CDR(\Lambda,X)$ and

\begin{enumerate}
\item[(ER1)] each function $y_i : [1,|y_i|] \to X^\pm$ is computable, that is, one can effectively
determine $y_i(\alpha)$ for every $\alpha \in [1,|y_i|]$ and $i \in [1,m]$,

\item[(ER2)] for every $i,j \in [1,m]$ and every $\alpha_i \in [1, |y_i|],\ \alpha_j \in [1, |y_j|]$
one can effectively compute $c(h_i,h_j)$, where $h_i = y_i^{\pm 1} \mid_{[\alpha_i,|y_i|]},\ h_j =
y_j^{\pm 1} \mid_{[\alpha_j,|y_j|]}$.
\end{enumerate}

Observe that since every $y_i$ is computable, $y_i^{-1}$ is computable too for every $i \in [1,m]$.
Next, it is obvious that concatenation of computable functions is computable, as well as restriction
of a computable function to a computable domain. Thus, if $g_i \ast g_j = h_i \circ h_j$, where
$g_i = y_i^{\delta_i} = h_i \circ c,\ g_j = y_j^{\delta_j} = c^{-1} \circ h_j,\ \delta_i, \delta_j
= \pm 1$, then both $h_i$ and $h_j$ are computable as restrictions $h_i = g_i \mid_{[1,\alpha]},\
h_j = g_j \mid_{[\alpha + 1,|g_j|]}$ for $\alpha = |c| = c(g_i^{-1},g_j)$, and so is $g_i \ast g_j$.
Now, using (ER2) twice we can determine $c((g_i \ast g_j)^{-1}, g_k)$, where $g_k = y_k^{\delta_k},
\ \delta_k = \pm 1$. Indeed, $c((g_i \ast g_j)^{-1}, g_k) = c(h_j^{-1} \circ h_i^{-1}, g_k)$, so,
if $c(h_j^{-1}, g_k) < |h_j^{-1}|$ then $c((g_i \ast g_j)^{-1}, g_k) = c(h_j^{-1}, g_k)$ which is
computable by (ER2), and if $c(h_j^{-1}, g_k) \geqslant |h_j^{-1}|$ then $c((g_i \ast g_j)^{-1},
g_k) = |h_j| + c(h_i^{-1}, h_k)$, where $h_k = g_k \mid_{[|h_j|+1, |g_k|]}$ -- again, all components
are computable and so is $c((g_i \ast g_j)^{-1}, g_k)$. It follows that $y_i^{\pm 1} \ast y_j^{\pm
1} \ast y_k^{\pm 1}$ is a computable function for every $i,j,k \in [1,m]$. Continuing in the same
way by induction one can show that every finite product of elements from $Y^{\pm 1}$, that is, every
element of $G$ given as a finite product of generators and their inverses, is computable as a
function defined over a computable segment in $\Lambda$ to $X^\pm$. Moreover, for any $g,h \in G$
one can effectively find $com(g,h)$ as a computable function. In particular, we automatically get a
solution to the Word Problem in $G$ provided $G$ has an effective representation by $\Lambda$-words
over an alphabet $X$.

\subsection{Computable universal trees}
\label{subs:comp_unive_trees}

Suppose $G$ has an effective representation by $\Lambda$-words over an alphabet $X$ and let $\Gamma_G$
be the universal $\Lambda$-tree for $G$. According to the construction in Section \ref{sec:universal},
every point of $\Gamma_G$ can be viewed as a pair $(\alpha, g)$, where $g \in G$ and $\alpha \in [0,
|g|]$. Such a pair is not unique but given another pair $(\beta, f)$ we can effectively find out if
both pair define the same point of $\Gamma_G$. Indeed, $(\alpha, g) \sim (\beta, f)$ if and only if
$\alpha = \beta \in [0, c(f,g)]$ and $c(f,g)$ can be found effectively.

Next, according to the definition given in Section \ref{sec:universal}, for $f \in G$ and $(\alpha,
g)$ representing a point in $\Gamma_G$, the image of $(\alpha, g)$ is defined as follows
$$f \cdot (\alpha, g) = (|f| + \alpha - 2 \min\{\alpha, c(f^{-1}, g)\}, f)$$
if $\alpha \leqslant c(f^{-1},g)$, and
$$f \cdot (\alpha, g) = (|f| + \alpha - 2 \min\{\alpha, c(f^{-1}, g)\}, f \ast g)$$
if $\alpha > c(f^{-1}, g)$. Since $G$ has an effective representation by $\Lambda$-words, it follows
that $c(f^{-1},g)$ can be found effectively and $f \ast g$ is a computable function. Thus, $f \cdot
(\alpha, g)$ can be determined effectively.

Summarizing the discussion above we prove the following result.

\begin{theorem}
\label{th:effect}
Let $G$ be a finitely generated group which has an effective representation by $\Lambda$-words
over an alphabet $X$ and let $\Gamma_G$ be the universal $\Lambda$-tree for $G$. Then
\begin{enumerate}
\item[(1)] $\Gamma_G$ can be constructed effectively in the sense that there exists a procedure
(infinite in general) which builds $\Gamma_G$ and every point of $\Gamma_G$ appears in the process
after finitely many steps,
\item[(2)] the action of $G$ on $\Gamma_G$ is effective in the sense that given a representative of
a point $v \in \Gamma_G$ and an element $g \in G$, one can effectively compute a representative of
$g \cdot v$.
\end{enumerate}
\end{theorem}
\begin{proof} Let $G = \langle Y \rangle$, where $Y$ is finite. The procedure building $\Gamma_G$
enumerates all finite words in the alphabet $Y^{\pm 1}$ and finds the corresponding computable
functions. According to the construction of Section \ref{sec:universal}, every point of $\Gamma_G$
can be viewed as a pair $(\alpha, g)$, where $g \in G$ and $\alpha \in [0, |g|]$, so eventually
every point of $\Gamma_G$ appears in the process. This concludes the proof of (1).

\smallskip

Finally, (2) follows from the discussion preceding the theorem.
\end{proof}

\bibliography{./../main_bibliography}

\end{document}